%
%
\input amstex
\documentstyle{tmna}
\NoBlackBoxes
\leftheadtext{Oleg Makarenkov} \rightheadtext{Perturbed
Poincar\'e-Andronov operators}

\topmatter
\title Influence of a small perturbation on Poincar\'e-Andronov operators with not well defined topological
degree
\endtitle
\author Oleg Makarenkov$^1$ \endauthor
\address {\smc Oleg Makarenkov}\\
Research Institute of Mathematics, Voronezh State University,
394006, Voronezh, Universitetskaja pl.1\endaddress
\curraddr Research Institute of Mathematics, Voronezh State
University,
394006, Voronezh, Universitetskaja pl.1\endcurraddr 
\email omakarenkov\@math.vsu.ru\endemail

\cvol{00} \cvolyear{0000} \cyear{0000} \mrecvd{October 10, 2007}

\submittedby J.~Mawhin\endsubmittedby

\subjclass Primary 54C40, 14E20; Secondary 46E25, 20C20\endsubjclass

\keywords Topological degree, perturbed Poincar\'e-Andronov map,
zero measure singularities\endkeywords

\abstract Let ${\Cal P}_\varepsilon\in C^0({\bold R}^n,{\bold
R}^n)$ be the Poincar\'e-Andronov operator over period $T>0$ of
 $T$-periodically perturbed autonomous system $\dot
x=f(x)+\varepsilon g(t,x,\varepsilon),$ where $\varepsilon>0$ is
small. Assuming that for $\varepsilon=0$ this system has a
$T$-periodic limit cycle $x_0$ we evaluate the topological degree
$d(I-{\Cal P}_\varepsilon,U)$ of $I-{\Cal P}_\varepsilon$ on an
open bounded set $U$ whose boundary $\partial U$ contains
$x_0([0,T])$ and ${\Cal P}_0(v)\not=v$ for any $v\in
\partial U\backslash x_0([0,T]).$  We give an explicit formula
connecting $d(I-{\Cal P}_\varepsilon,U)$ with the topological
indices of zeros of the associated Malkin's bifurcation function.
The goal of the paper is to prove the Mawhin's conjecture claiming
that $d(I-{\Cal P}_\varepsilon,U)$ can be any integer in spite of
the fact that the measure of the set of fixed points of ${\Cal
P}_0$ on $\partial U$ is zero.\endabstract

\thanks The work is partially supported by the Grant BF6M10 of Russian
Federation Ministry of Education and CRDF (US), and by RFBR Grant
06-01-72552, and by the President of Russian Federation Young
Candidate of Science grant MK-1620.2008.1.
\endthanks
\endtopmatter

\document

\head 1. Introduction\endhead

Consider the system of ordinary differential equations
$$
   \dot x=f(x)+\varepsilon g(t,x,\varepsilon),
   \tag 1.1
$$
where $f\in C^1({\bold R}^n,{\bold R}^n),$ $g\in C^0({\bold
R}\times{\bold R}^n\times[0,1],{\bold R}^n),$
$g(t+T,v,\varepsilon)\equiv g(t,v,\varepsilon)$ and
$\varepsilon>0$ is a small parameter. We suppose that equation
(1.1) defines a flow in ${\bold R}^n,$ i.e. assume the uniqueness
and global existence for the solutions of the Cauchy problems
associated to (1.1). For each $v\in{\bold R}^n$ we denote by
$x_\varepsilon(\cdot,v)$ the solution of (1.1) with
$x_\varepsilon(0,v)=v.$ Thus, the Poincar\'e-Andronov operator
over the period $T>0$ is defined by
$$
  {\Cal P}_\varepsilon(v):=x_\varepsilon(T,v).
$$
The problem of the existence
 (and even stability, see Ortega
\cite{11}) of $T$-periodic solutions of (1.1) with initial
conditions inside an open bounded set $U$ can be solved by
evaluating the topological degree $d(I-{\Cal P}_\varepsilon,U)$ of
$I-{\Cal P}_\varepsilon$ on $U$ (see \cite{6}). In the case when
${\Cal P}_0$ has no fixed points on the boundary $\partial U$ of
$U$ the problem is completely solved by Capietto, Mawhin and
Zanolin \cite{2} who proved that $d(I-{\Cal P}_0,U)=(-1)^n d(f,U)$
generalizing the result by Berstein and Halanay \cite{1} where $U$
is assumed to be a neighborhood of an isolated zero of $f.$ In the
case when ${\Cal P}_0$ has fixed points on $\partial U$ the
pioneer result has been obtained by Mawhin \cite{10} who
considered the situation when $f=0.$ Mawhin proved that if
$g_0(v)=\int_0^T g(\tau,v,0)d\tau$ does not vanish on $\partial U$
then $d(I-{\Cal P}_\varepsilon,U)$ is defined for $\varepsilon>0$
sufficiently small and it can be evaluated as $d(I-{\Cal
P}_\varepsilon,U)=d(-g_0,U).$ This paper studies an intermediate
situation when the fixed points of ${\Cal P}_0$ fill a part of
$\partial U.$ Current results on this subject deal with the case
when $\partial U$ contains a fixed number of fixed points, e.g.
Feckan \cite{4}, Kamenskii-Makarenkov-Nistri \cite{5}.
 As a part of a wider study of
this problem Jean Mawhin (his seminar, November~2005) asked a
question on evaluating  $d(I-{\Cal P}_\varepsilon,U)$ in the case
when $\partial U$ contains a curve of fixed points of ${\Cal
P}_0.$ He settled  the following conjecture:

\vskip0.1cm

\noindent {\bf Mawhin's conjecture.} {\it For small
$\varepsilon>0$ the topological degree $d(I-{\Cal
P}_\varepsilon,U)$ can be any integer depending on the
perturbation term $g$ in spite of the fact that the measure of
$\{v\in\partial U:{\Cal P}_0(v)=v\}$ is zero.}

\vskip0.1cm

 The goal of this paper is to evaluate $d(I-{\Cal P}_\varepsilon,U)$ and to
 give a proof of the above conjecture in the case when $\{v\in\partial U:{\Cal P}_0(v)=v\}$
 forms a curve coming from a $T$-periodic limit cycle of the unperturbed system
$$
   \dot x=f(x).
   \tag 1.2
$$
Our fundamental assumption is that the algebraic multiplicity of
the multiplicator $+1$ of the linearized system
$$
  \dot y=f'(x_0(t))y
  \tag 1.3
$$
equals to $1.$ In this case we say that the cycle $x_0$ is
nondegenerate.

The paper is organized as follows. In Section~2 for a fixed point
$v_\varepsilon$ of ${\Cal P}_\varepsilon$ satisfying
$v_\varepsilon\to v_0\in x_0([0,T])$ as $\varepsilon\to 0$ we
obtain an asymptotic direction of the vector $v_\varepsilon-v_0.$
By means of this result we evaluate in Section~3 the topological
index of such fixed points $v_\varepsilon\to v_0\in x_0([0,T])$ as
$\varepsilon\to 0$ that $v_\varepsilon\in U.$  Finally in
Section~4 we give a proof of the Mawhin's conjecture provided that
a technical assumption (see assumption~4.1) is satisfied.

\head 2. Direction the fixed points of Poincar\'e-Andronov
operator move when the perturbation increases\endhead

Since the cycle $x_0$ is nondegenerate we can define (see
\cite{3}, Ch.~IV, \S~20, Lemma~1)  a matrix function $Z_{n-1}$
solving the adjoint system
$$
  \dot z=-(f'(x_0(t)))^*z
  \tag 2.1
$$
and having the form $Z_{n-1}(t)=\Phi(t){\roman e}^{\Lambda t},$
where $\Phi$ is a continuous $T$-periodic $n\times n-1$ matrix
function and $\Lambda$ is a $n-1\times n-1$-matrix with different
from $0$ eigenvalues. Let  $z_0$ be the $T$-periodic solution of
(2.1) satisfying $z_0(0)^*\dot x_0(0)=1.$ Finally, we denote by
$Y_{n-1}$ the $n\times n-1$ matrix function whose columns are
solutions of the linearized system (1.3) satisfying
$Y_{n-1}(0)^*Z_{n-1}(0)=I.$

\vskip0.3cm

 The results of this paper are formulated in terms of the following auxiliary functions:
$$   M(\theta) = \int\limits_{0}^{T}
   z_0(\tau)^*g(\tau-\theta,{x}_0(\tau),0)d\tau,$$ $$
    M^\bot(t,\theta)= \left({\roman e}^{\Lambda T}\right)^*\left(\left({\roman e}^{\Lambda
T}\right)^*-I\right)^{-1}\int\limits_{t-T+\theta}^{t+\theta}
   \left({Z_{n-1}(\tau)}\right)^*g(\tau-\theta,{x}_0(\tau),0)d\tau,$$
   $$
  \angle(u,v)=\arccos\frac{\left<u,v\right>}{\|u\|\cdot\|v\|}.
$$

The function $M$ was proposed by Malkin (see \cite{9},
formula~3.13) and the function $M^\bot$ is a generalization of the
function $M^\bot_z$ of \cite{8}.

Next Theorem~2.1 shows that if a family
$\{x_{\varepsilon,\lambda}\}_{\lambda\in\Lambda}$ of $T$-periodic
solutions of (1.1) emanate from $x_0(\cdot+\theta_0)$ then a
suitable projection of
$x_{\varepsilon,\lambda}(t)-{x}_0(t+\theta_0)$ can be always
controlled. Though motivated by the Mawhin's conjecture,
Theorem~2.1
can be of a general interest in the theory of oscillations playing
a role of the first approximation formula  (see Loud \cite{7},
formula~1.3, Lemma~1 and formula for $x$ at p.~510) in the case
when the zeros of the  bifurcation function $M$ are not necessary
isolated.

\proclaim{Theorem 2.1} Let $x_0$ be a nondegenerate $T$-periodic
cycle of (1.2). Let
$\{x_{\varepsilon,\lambda}\}_{\lambda\in\Lambda}$ be a family of
$T$-periodic solutions of (1.1) such that
$x_{\varepsilon,\lambda}(t)\to x_0(t+\theta_0)$ as $\varepsilon\to
0$ uniformly with respect to $t\in[0,T]$ and $\lambda\in\Lambda.$
Then
$$
  \angle\left(Z_{n-1}(t+\theta_0)^*(
  x_{\varepsilon,\lambda}(t)-{x}_0(t+\theta_0)),M^\bot(t,\theta_0)\right)\to 0\quad{as\
  }\varepsilon\to 0
$$ uniformly with respect to $t\in[0,T]$ and $\lambda\in\Lambda.$
\endproclaim

\noindent{\bf Proof.} The proof makes use of the idea of
Theorem~3.1 of \cite{8}. In the sequel $(A,B)$ denotes the matrix
composed by columns of matrixes $A$ and $B.$
 Let
  $a_\varepsilon\in C^0([0,T],{\bold R}^n)$ be given by
$$
  a_{\varepsilon}(t)=(z_0(t+\theta_0),Z_{n-1}(t+\theta_0))^*\,(x_\varepsilon(t)-{x}_0(t+\theta_0)).
  \tag 2.2
$$
Denoting $Y(t)=(\dot x_0(t),Y_{n-1}(t))$ by Perron's lemma
\cite{12} (see also Demidovich (\cite{3}, Sec. III, \S 12) we have
$$
  (z_0(t),Z_{n-1}(t))^*\, Y(t)=I,\quad{\roman f\roman o\roman r\ \roman a\roman n\roman y\
  }t\in{\bold R}.
$$
Thus
$$
x_\varepsilon(t)-{x}_0(t+\theta_0)=Y(t+\theta_0)a_\varepsilon(t),\quad{\roman
f\roman o\roman r\ \roman a\roman n\roman y\
  }t\in{\bold R}.
  \tag 2.3
$$

 By subtracting (1.2) where $x$ is replaced by $x_0(\cdot+\theta_0)$  from (1.1) where $x$ is replaced by
$x_\varepsilon$ we obtain
$$
  \dot x_\varepsilon(t)-\dot{
  {x}}_0(t+\theta_0)=f'({x}_0(t+\theta_0))(x_\varepsilon(t)-{x}_0(t+\theta_0))$$
  $$
    +\varepsilon g(t,x_\varepsilon(t),\varepsilon)+o(t,x_\varepsilon(t)-{x}_0(t+\theta_0)),
\tag 2.4
$$
where $o(t,v)/\|v\|\to 0$ as ${\bold R}^n\ni v\to 0$ uniformly
with respect to $t\in[0,T].$ By substituting (2.3) into (2.4) we
have
$$
    \dot Y(t+\theta_0)a_\varepsilon(t)+  Y(t+\theta_0)\dot a_\varepsilon(t)
$$
$$
 =
 f'({x}_0(t+\theta_0))Y(t+\theta_0)a_\varepsilon(t)+\varepsilon g(t,x_\varepsilon(t),\varepsilon)+
  o(t,x_\varepsilon(t)-{x}_0(t+\theta_0)).
$$
Since $f'({x}_0(t))Y(t)=\dot Y(t)$ the last relation can be
rewritten as
$$
  Y(t+\theta_0)\dot a_\varepsilon(t)
 =\varepsilon g(t,x_\varepsilon(t),\varepsilon)+
  o(t,x_\varepsilon(t)-{x}_0(t+\theta_0)).
 \tag 2.5
$$
Applying $Z_{n-1}(t+\theta_0)^*$ to both sides of (2.5) we have
$$
  (0,I)\dot a_{\varepsilon}(t)=\varepsilon Z_{n-1}(t+\theta_0)^*\, g(t,x_\varepsilon(t),\varepsilon)
  +Z_{n-1}(t+\theta_0)^*\, o(t,x_\varepsilon(t)-{x}_0(t+\theta_0)),
$$
where $0$ denotes the $n-1$ dimensional zero vector  and $I$ stays
for the identical $n-1\times n-1$ matrix. So
$$
  (0,I)a_{\varepsilon}(t)=(0,I)a_{\varepsilon}(t_0)+\varepsilon\int\limits_{t_0}^t
  Z_{n-1}(\tau+\theta_0)^*g(\tau,
  x_\varepsilon(\tau),\varepsilon)d\tau
$$ $$   +\int\limits_{t_0}^t
Z_{n-1}(\tau+\theta_0)^*o(\tau,x_\varepsilon(\tau)-x_0(\tau+\theta_0))d\tau.
\tag 2.6 $$ From the definition of $Z_{n-1}$ we have that
$Z_{n-1}(t)^*=\left({\roman e}^{\Lambda T}\right)^*Z_{n-1}(t-T)^*$
for any $t\in{\bold R}$ and so $(0,I)a_{\varepsilon}(t)$ satisfies
$$
  (0,I)a_{\varepsilon}(t_0)= \left({\roman e}^{\Lambda
T}\right)^* (0,I) a_{\varepsilon}(t_0-T)\quad{\roman f\roman
o\roman r\ \roman a\roman n\roman y\
  }t_0\in[0,T].
  \tag 2.7
$$
Solving (2.6)-(2.7) with respect to $(0,I)a_{\varepsilon,n}(t_0)$
we obtain
$$
  (0,I)a_{\varepsilon}(t_0)
=\varepsilon  \left({\roman e}^{\Lambda
T}\right)^*\left(\left({\roman e}^{\Lambda
T}\right)^*-I\right)^{-1} \int\limits_{t_0-T}^{t_0}
Z_{n-1}(\tau+\theta_0)^*g(\tau,
  x_\varepsilon(\tau),\varepsilon)d\tau $$ $$
   \hskip-1cm+\left({\roman e}^{\Lambda T}\right)^*\left(\left({\roman e}^{\Lambda T}\right)^*-I\right)^{-1}\int\limits_{t_0-T}^{t_0}
  Z_{n-1}(\tau+\theta_0)^*{o(\tau,x_\varepsilon(\tau)-x_0(\tau+\theta_0))}d\tau
$$
for any $t_0\in[0,T].$ On the other hand from (2.2) we obtain $$
Z_{n-1}(t+\theta_0)^*( x_\varepsilon(t)-{x}_0(t+\theta_0))= (0,I)
a_{\varepsilon}(t)$$ and therefore
$$
  Z_{n-1}(t+\theta_0)^*( x_\varepsilon(t)-{x}_0(t+\theta_0))
  -q_\varepsilon(t)$$
$$
   =\varepsilon  \left({\roman e}^{\Lambda T}\right)^*\left(\left({\roman e}^{\Lambda T}\right)^*-I\right)^{-1} \int\limits_{t-T}^{t}
Z_{n-1}(\tau+\theta_0)^*g(\tau,
  x_\varepsilon(\tau),\varepsilon)d\tau,
  \tag 2.8
$$
where $$
 q_\varepsilon=\left({\roman e}^{\Lambda T}\right)^*\left(\left({\roman e}^{\Lambda T}\right)^*-I\right)^{-1}\int\limits_{t-T}^{t}
  Z_{n-1}(\tau+\theta_0)^*{o(\tau,x_\varepsilon(\tau)-x_0(\tau+\theta_0))}d\tau.
$$
From (2.8) we obtain
$$
 \angle\left(Z_{n-1}(t+\theta_0)^*(
 x_\varepsilon(t)-{x}_0(t+\theta_0)),M^\bot(t,\theta_0)\right)$$
 $$
  = \angle\left(Z_{n-1}(t+\theta_0)^*\dfrac{x_\varepsilon(t)-{x}_0(t+\theta_0)}
 {\|x_\varepsilon-x_0(\cdot+\theta_0)\|_{[0,T]}},M^\bot(t,\theta_0)\right)$$
 $$
  - \angle\left(Z_{n-1}(t+\theta_0)^*\dfrac{x_\varepsilon(t)-{x}_0(t+\theta_0)}
 {\|x_\varepsilon-x_0(\cdot+\theta_0)\|_{[0,T]}}-\dfrac{q_\varepsilon(t)}{\|x_\varepsilon-x_0(\cdot+\theta_0)\|_{[0,T]}},M^\bot(t,\theta_0)\right)$$
 $$
 +
 \angle\left(\left({\roman e}^{\Lambda T}\right)^*\left(\left({\roman e}^{\Lambda T}\right)^*-I\right)^{-1} \int\limits_{t-T}^{t}
Z_{n-1}(\tau+\theta_0)^*g(\tau,
  x_\varepsilon(\tau),\varepsilon)d\tau,M^\bot(t,\theta_0)\right).
$$
But the difference of the first two terms in the right hand part
of the last equality tends to zero as $\varepsilon\to 0$ and thus
the thesis follows.

\qed

\

\noindent Next theorem is a reformulation of  Theorem~2.1 suitable
for our further considerations.

\proclaim{Theorem 2.2} Let $x_0$ be a nondegenerate $T$-periodic
cycle of (1.2). Let
$\{x_{\varepsilon,\lambda}\}_{\lambda\in\Lambda}$ be a family of
$T$-periodic solutions of (1.1) such that
$x_{\varepsilon,\lambda}(t)\to x_0(t+\theta_0)$ as $\varepsilon\to
0$ uniformly with respect to $t\in[0,T]$ and $\lambda\in\Lambda.$
Let $l\in{\bold R}^n$ be an arbitrary vector such that
$\left<l,\dot x_0(\theta_0)\right>=0.$ Assume that
$\left<l,Y_{n-1}(\theta_0)M^\bot(0,\theta_0)\right>\not=0.$ Then
there exists $\varepsilon_0>0$ such that
$$\left<l,x_{\varepsilon,\lambda}(0)-x_0(\theta_0)\right>>0\quad or \quad
\left<l,x_{\varepsilon,\lambda}(0)-x_0(\theta_0)\right><0$$
according as
$$\left<l,Y_{n-1}(\theta_0)M^\bot(0,\theta_0)\right>>0\quad or
\quad\left<l,Y_{n-1}(\theta_0)M^\bot(0,\theta_0)\right><0$$ for
any $\lambda\in\Lambda$ and any $\varepsilon\in(0,\varepsilon_0].$
\endproclaim

\noindent{\bf Proof.}  By Perron's lemma \cite{12} (see also
Demidovich (\cite{3}, Sec. III, \S 12) we have
  $$
    v=Y_{n-1}(t)Z_{n-1}(t)^*v+\dot x_0(t)z_0(t)^*v
  $$
for any $v\in{\bold R}^n$ and $t\in{\bold R}.$ Therefore

\centerline{$
\left<l,x_{\varepsilon,\lambda}(0)-x_0(\theta_0)\right>$}

\centerline{$
       =\left<l,Y_{n-1}(\theta_0)Z_{n-1}(\theta_0)^*(x_{\varepsilon,\lambda}(0)-x_0(\theta_0))\right.$}

   \centerline{$+\left.\dot
x_0(\theta_0)z_0(\theta_0)^*(x_{\varepsilon,\lambda}(0)-x_0(\theta_0))
\right>$}

\centerline{$\left<Y_{n-1}(\theta_0)^*l,Z_{n-1}(\theta_0)^*(x_{\varepsilon,\lambda}(0)-x_0(\theta_0))\right>.$}

Since $\left<Y_{n-1}(\theta_0)^*l,M^\bot(0,\theta_0)\right>\not=0$
then by Theorem~2.1 there exists $\varepsilon_0>0$ such that
$${\roman s\roman i\roman g\roman n}
\left<Y_{n-1}(\theta_0)^*l,Z_{n-1}(\theta_0)^*(x_{\varepsilon,\lambda}(0)-x_0(\theta_0))\right>={\roman
s\roman i\roman g\roman
n}\left<Y_{n-1}(\theta_0)^*l,M^\bot(0,\theta_0)\right>$$ for any
$\varepsilon\in(0,\varepsilon_0]$ and $\lambda\in\Lambda$ and thus
the proof is complete. \qed

\head 3. The topological degree of the perturbed
Poincar\'e-Andronov operator \endhead

To proceed to the proof of our main Theorem~3.1 we need three
additional theorems which are formulated below for the convenience
of the reader.

 {\bf Malkin's Theorem} (see \cite{9}, p.~41) {\it Assume that $T$-periodic solutions $x_\varepsilon$ of (1.1)
 satisfy
 the property $x_\varepsilon(t)\to
x_0(t+\theta_0)$ as $\varepsilon\to 0.$ Then $M(\theta_0)=0.$}

 {\bf
Capietto-Mawhin-Zanolin Theorem} (see \cite{2}, Corollary~2). {\it
Let $V\subset{\bold R}^n$ be an open bounded set. Assume that
${\Cal P}_0(v)\not=v$ for any $v\in\partial V.$ Then $d(I-{\Cal
P}_0,V)=(-1)^n d(f,V).$}

{\bf Kamenskii-Makarenkov-Nistri Theorem} (see \cite{5},
Corollary~2.8). {\it Assume that $\theta_0\in[0,T]$ is an isolated
zero of the bifurcation function $M.$ Then there exist
$\varepsilon_0>0$ and  $r>0$ such that ${\Cal
P}_\varepsilon(v)\not=v$ for any $\|v-v_0\|=r$ and any
$\varepsilon\in(0,\varepsilon_0].$ Moreover $d(I-{\Cal
P}_\varepsilon,B_r(v_0))={\roman i\roman n\roman d}(\theta_0,M).$}

\

We will say that the  set $U\subset{\bold R}^n$ has a smooth
boundary if given any $v\in\partial U$ there exists $r>0$ and a
homeomorphism of $\{\xi\in{\bold R}^{n-1}:\|\xi\|\le 1\}$ onto
$\partial U\cap B_r(v).$ Thus any set $U$ with a smooth boundary
possesses a tangent plane to $\partial U$ at any $v\in\partial U.$
This tangent plane will be denoted by $L_U(v).$ Moreover, if $U$
has a smooth boundary and ${\bold R}^n\ni h\not\in L_U(v)$ then
there exists $\lambda_0>0$ such that either $\lambda h+v\in U$ for
any $\lambda\in(0,\lambda_0]$ or $\lambda h+v\not\in U$ for any
$\lambda\in(0,\lambda_0].$ In this case we will say that $h$
centered at $v$ is directed inward to $U$ or outward respectively.

\proclaim{Theorem 3.1}  Let $x_0$ be a nondegenerate $T$-periodic
cycle of (1.2). Let $U\subset{\bold R}^n$ be an open bounded set
with a smooth boundary and $x_0([0,T])\subset\partial U.$ Assume
that ${\Cal P}_0(v)\not=v$ for any $v\in\partial U\backslash
x_0([0,T]).$ Assume that $M$ has a finite number of zeros $0\le
\theta_1<\theta_2<...<\theta_k<T$ on $[0,T]$ and ${\roman i\roman
n\roman d}(\theta_i,M)\not=0$ for any $i\in\overline{1,k}.$ Assume
that $Y_{n-1}(\theta_i)M^\bot(0,\theta_i)\not\in
L_U(x_0(\theta_i))$ for any $i\in\overline{1,k}.$ Then for any
$\varepsilon>0$ sufficiently small $ d(I-{\Cal P}_\varepsilon,U)$
is defined. Moreover,
$$
  d(I-{\Cal P}_\varepsilon,U)=(-1)^n d(f,U)-\sum_{i=1}^k {\roman i\roman n\roman d}(\theta_i,M)D_i,
$$
where $D_i=1$ or $D_i=0$  according as
$Y_{n-1}(\theta_i)M^\bot(0,\theta_i)$ centered at $x_0(\theta_i)$
is directed inward to $U$ or outward.
\endproclaim

\noindent{\bf Proof.} By Kamenskii-Makarenkov-Nistri theorem there
exists $r>0$ and $\varepsilon_0>0$ such that
$$ 
  d(I-{\Cal P}_\varepsilon,B_r(x_0(\theta_i)))={\roman i\roman n\roman d}(\theta_i,M)
  \tag 3.1
$$
  for any $\varepsilon\in(0,\varepsilon_0]$ and $i\in\overline{1,k}.$
From Malkin's theorem we have the following "Malkin's property":
$r>0$ can be decreased, if necessary, in such a way that there
exists $\varepsilon_0>0$ such that any $T$-periodic solution
$x_\varepsilon$ of (1.1) with initial condition
$x_\varepsilon(0)\in B_r(x_0([0,T]))$ and
$\varepsilon\in(0,\varepsilon_0]$ satisfies $x_\varepsilon(0)\in
\cup_{i\in\overline{1,k}} B_r (x_0(\theta_i)).$ Malkin's property
implies that
$$
  d\left(I-{\Cal P}_\varepsilon,\left(B_r(x_0([0,T]))\backslash\cup_{i\in\overline{1,k}}B_r(x_0(\theta_i))\right)\cap U\right)=0
  \tag 3.2
$$
  for any $\varepsilon\in(0,\varepsilon_0].$
 Denote by $l_i$ the perpendicular to $L_U(x_0(\theta_i))$
directed outward away from $U$ or inward according as
$(Z_{n-1}(\theta_i)^*)^{-1}M^\bot(0,\theta_i)$ centered at
$x_0(\theta_i)$ is directed outward away from $U$ or inward. From
Theorem~2.2 and Malkin's property we have that $\varepsilon_0>0$
can be diminished in such a way that for any $i\in\overline{1,k}$
any $T$-periodic solution $x_\varepsilon$ of (1.1) with initial
condition $x_\varepsilon(0)\in B_r(x_0(\theta_i))$ and
$\varepsilon\in(0,\varepsilon_0]$ satisfies $x_\varepsilon(0)\in
B_r (x_0(\theta_i))\cap U$ or $x_\varepsilon(0)\not\in B_r
(x_0(\theta_i))\cap U$ according as $D_i=1$ or $D_i=0.$ This
observation allows to deduce from (3.1) that
$$
  { d}(I-{\Cal P}_\varepsilon,B_r (x_0(\theta_i))\cap
  U)={\roman i\roman n\roman d}(\theta_i,M),\ \ {\roman i\roman f}\ D(\theta_i)=1,
  \tag 3.3
$$
$$ 
  { d}(I-{\Cal P}_\varepsilon,B_r (x_0(\theta_i))\cap
  U)=
    0,\ \ {\roman i\roman f}\ D(\theta_i)=0,
  \tag 3.4
$$

\noindent for any $\varepsilon\in(0,\varepsilon_0]$ and
$i\in\overline{1,k}.$

Observe that our choice of $r>0$ ensures that ${\Cal
P}_0(v)\not=v$ for any $v\in\partial\left(U\backslash
B_r(x_0([0,T]))\right).$ Thus, by Capietto-Mawhin-Zanolin theorem
we have $
  d(I-{\Cal P}_0,U\backslash B_r(x_0([0,T])))=(-1)^n d(f,U\backslash
  B_r(x_0([0,T]))).
$ Without loss of generality we can consider $r>0$ sufficiently
small such that $d(f,U\backslash
  B_r(x_0([0,T])))=d(f,U)$ obtaining
$$
  d(I-{\Cal P}_0,U\backslash B_r(x_0([0,T])))=(-1)^n d(f,U).
\tag 3.5
$$
Since
$$
  d(I-{\Cal P}_\varepsilon,U)=d\left(I-{\Cal
P}_\varepsilon,\left(B_r(x_0([0,T]))\backslash\cup_{i\in\overline{1,k}}B_r(x_0(\theta_i))\right)\cap
  U\right)
$$
$$
  + { d}\left(I-{\Cal P}_\varepsilon,\cup_{i\in\overline{1,k}}B_r (x_0(\theta_i))\cap
  U\right)
$$
$$
   + d(I-{\Cal P}_\varepsilon,U\backslash B_r(x_0([0,T])))
$$
the conclusion follows from formulas (3.2)-(3.5). \qed

\head 4. A proof of the Mawhin's conjecture \endhead

In this section we assume that the set $U\subset{\bold R^n}$ has a
smooth boundary and there exists $v_{n-1}\in\bold R^{n-1}$
satisfying the following assumption
$$
Y_{n-1}(t)\left({\roman e}^{\Lambda
     T}\right)^*\left(\left({\roman
     e}^{\Lambda T}\right)^*-I\right)^{-1}\left({\roman e}^{\Lambda
     t}\right)^*v_{n-1}\not\in L_U(t)\quad{\roman f\roman o\roman r \ \roman a\roman n\roman y\ } t\in[0,T]. \tag 4.1
$$
We note that assumption (4.1) does not depend on the perturbation
term of (1.1) and relies to unperturbed system (1.2). Let $D=1$ or
$D=0$ according as $Y_{n-1}(0)\left({\roman e}^{\Lambda
     T}\right)^*\left(\left({\roman
     e}^{\Lambda T}\right)^*-I\right)^{-1}\left({\roman e}^{\Lambda
     t}\right)^*v_{n-1}$ centered at $x_0(0)$ is directed inward to $U$ or
     outward.
Given odd $m\in\bold N$
 we construct the perturbation term $g$
in such a way that $d(I-{\Cal P}_\varepsilon,U)=(-1)^n d(f,U)
-m(2D-1)$ for any $\varepsilon>0$ sufficiently small. Without loss
of generality we consider $T=2\pi.$

Since $(z_0(t),Z_{n-1}(t))$ is nonsingular then
$((z_0(t),\Phi(t))^*$ is nonsingular as well. Define
$\Omega:x_0([0,2\pi])\to\bold R^n$ as
$\Omega(x_0(t))=((z_0(t),\Phi(t))^*)^{-1}$ for any $t\in[0,2\pi].$
By Uryson's theorem (see \cite{6}, Ch.~1,~Theorem~1.1) $\Omega$
can be continued to the whole $\bold R^n$ in such a way that
$\Omega\in C^0(\bold R^n,\bold R^n).$ Analogously, we consider
$\widetilde{\Gamma}\in C^0(\bold R^n,\bold R^n)$ such that
$\widetilde{\Gamma}(x_0(t))=\left(\arcsin(\sin t),0,\hdots,
0\right)^*$ and denote by $\Gamma\in C^0(\bold R^n,\bold R)$ the
first component of $\widetilde{\Gamma}.$ Let us define a
$2\pi$-periodic $\alpha$-approximation of $\left(\left({\roman
e}^{\Lambda t}\right)^*\right)^{-1}$ on $[-2\pi,0]$ by
$$
  {\roman e}_\alpha(t)=(({\roman e}^{\Lambda t})^*)^{-1}, \ \  {\roman i\roman f\
    }t\in[-2\pi,-\alpha],
$$
$$
  {\roman e}_\alpha(t)=\frac{t}{-\alpha}\left(\left({\roman
    e}^{-\Lambda\alpha}\right)^*\right)^{-1}+\left(1-\frac{t}{-\alpha}\right)
    \left(\left({\roman
    e}^{-2\pi\Lambda}\right)^*\right)^{-1}, \ \ {\roman i\roman f\
    }t\in[-\alpha,0],
$$

\noindent which is continued to $(-\infty,\infty)$ by the
$2\pi$-periodicity.  We are now in a position to introduce the
required perturbation term, namely we consider that the perturbed
system (1.1) has the following form
$$
   \dot x =f(x)  +\varepsilon\Gamma(x)\Omega(x){
     D\sin(mt)+(1-D)\cos(mt)\choose
     (D\cos(mt)+(1-D)\sin(mt)){\roman e}_\alpha(t)v_{n-1}},
\tag 4.2
$$
where $\alpha>0$ is sufficiently small. Consequently we denote by
${\Cal P}_\varepsilon$ the Poincar\'e-Andronov operator of system
(4.2) over the period $2\pi.$

\proclaim{Proposition 4.1} Let $x_0([0,T])\subset U\subset\bold
R^n$ be an open bounded set with a smooth boundary and assume that
there exists $v_{n-1}\in\bold R^n$ such that (4.1) is satisfied.
 Then given any odd $m>0$
there exists $\alpha_0>0$ such that for any fixed
$\alpha\in(0,\alpha_0]$ and $\varepsilon>0$ sufficiently small
$d(I-{\Cal P}_\varepsilon,U)$ is defined and
$$
d(I-{\Cal P}_\varepsilon,U)=\left\{(-1)^n d(f,U)-m,\ \ {\roman
i\roman f\ }D=1, \atop (-1)^n d(f,U)+m, \ \  {\roman i \roman f\
}D=0. \right.
$$
\endproclaim

\noindent{\bf Proof.} By the definition of $\Omega$ and $\Gamma$
we have
$$
  {
    z_0(t)^*\choose
    Z_{n-1}(t)^*}\Omega(x_0(t))={1 \ \ \ \ \ \
    0 \ \ \choose 0 \ \ \left({\roman e}^{\Lambda
    t}\right)^*},
  \tag 4.3
$$
$$
    \Gamma(x_0(t))=\arcsin(\sin t).
$$
Therefore, taking into account that $m$ is odd, we obtain the
following formula for the bifurcation function $M$
$$
  M(\theta)=\int_0^{2\pi}\arcsin(\sin\tau)(D\sin(m(\tau-\theta))+(1-D)\cos(m(\tau-\theta)))d\tau\atop
  =(-1)^{(m-1)/2}\dfrac{4D\cos(m\theta)+4(1-D)\sin(m\theta)}{m^2}
$$
whose zeros are
$\theta_j=\frac{1}{m}\left(\frac{D\pi}{2}+j\pi\right),$
$j\in\overline{0,2m-1}.$ Moreover,
$$
{\roman i\roman n\roman d}(\theta_j,M)={\roman s\roman i\roman
g\roman n}(M'(\theta_j))
  \tag 4.4
$$
$$ =(-1)^{(m-1)/2}{\roman s\roman i\roman g\roman n}\left(\dfrac{4m\left(\hskip-0.1cm-D\sin\left(D\dfrac{\pi}{2}+j\pi\right)+(1\hskip-0.1cm-\hskip-0.1cmD)
  \cos\left(D\dfrac{\pi}{2}+j\pi\right)\right)}{m^2}\right).
$$
Let us denote by $M^\bot_\alpha$ the function $M^\bot$
corresponding to system (4.2). From (4.3) we have that
$$
M^\bot_\alpha(0,\theta)= \left({\roman e}^{\Lambda
     T}\right)^*\left(\left({\roman
     e}^{\Lambda T}\right)^*-I\right)^{-1}\int\limits_{-2\pi}^0(Z_{n-1}(s+\theta))^*g(s,x_0(s+\theta),0)ds
$$
$$
 =\left({\roman e}^{\Lambda
     T}\right)^*\left(\left({\roman
     e}^{\Lambda T}\right)^*-I\right)^{-1}\left({\roman e}^{\Lambda
     \theta}\right)^*\circ
$$
$$
      \circ\int\limits_{-2\pi}^0\left({\roman
     e}^{\Lambda s}\right)^*{\roman
     e}_\alpha(s)v_{n-1}\arcsin(\sin(s+\theta))(D\cos(ms)+(1-D)\sin(ms))ds.
$$
Since
$$
  \int\limits_{-2\pi}^0\arcsin(\sin(s+\theta))(D\cos(ms)+(1-D)\sin(ms))ds
$$
$$
  =
     -(-1)^{(m-1)/2}\cdot\dfrac{4(D\sin\left(m\theta\right)+(1-D)\cos(m\theta))}{m^2}
$$
by  taking into account that $m$ is odd we have that $
  M_\alpha^\bot(0,\theta)\to M^\bot_0(0,\theta)
$ as $\alpha\to 0,$ where
$$
M^\bot_0(0,\theta)=-\left({\roman e}^{\Lambda
     T}\right)^*\left(\left({\roman
     e}^{\Lambda T}\right)^*-I\right)^{-1}\left({\roman e}^{\Lambda
     \theta}\right)^*v_{n-1}(-1)^{(m-1)/2}
$$
$$
  \cdot\dfrac{4(D\sin\left(m\theta\right)+(1-D)\cos(m\theta))}{m^2}.
$$
Put
$q(\theta)=-(-1)^{(m-1)/2}(D\sin(m\theta)+(1-D)\cos(m\theta)).$
Then, taking any $\theta\in[0,2\pi]$ and using the definition of
$D$ we conclude that $Y_{n-1}(\theta)M^\bot_0(0,\theta)$ centered
at $x_0(\theta)$ is directed inward to $U$ or outward according as
${\roman s\roman i\roman g\roman n}(q(\theta))(2D-1)=1$ or
${\roman s\roman i\roman g\roman n}(q(\theta))(2D-1)=-1.$
Therefore, there exists $\alpha_0>0$ such that for any
$\alpha\in[0,\alpha_0]$ and any $\theta\in[0,2\pi]$ we have that
$Y_{n-1}(\theta)M^\bot_\alpha(0,\theta)$ centered at $x_0(\theta)$
is directed inward to $U$ or outward according as ${\roman s\roman
i\roman g\roman n}(q(\theta))(2D-1)=1$ or ${\roman s\roman i\roman
g\roman n}(q(\theta))(2D-1)=-1.$ Thus denoting by ${\Cal
P}_{\varepsilon,\alpha}$ the Poincar\'e-Andronov operator of
system (4.2) from Theorem~3.1 we have that
$$
  d(I-{\Cal P}_{\varepsilon,\alpha},U)=(-1)^n d(f,U)-\sum_{j\in\overline{0,2m-1}:{\roman s\roman i\roman g\roman n}(q(\theta_j))(2D-1)=1}{\roman i\roman n\roman d}(\theta_j,M)
  \tag 4.5
$$
for any $\alpha\in(0,\alpha_0].$ Consider the case when $D=1.$
Then the property \linebreak ${\roman s\roman i\roman g\roman
n}(q(\theta_j))(2D-1)=1$ is equivalent to
$$
(-1)^{(m-1)/2}{\roman s\roman i\roman g\roman
n}(\sin(\pi/2+j\pi))=-1. \tag 4.6
$$
If $j\in\overline{0,2m-1}$ satisfies (4.6)
  then
  (4.4) implies  ${\roman i\roman n\roman d}(\theta_j,M)=1.$ Since
  there exists exactly $m$ elements of
  $\overline{0,2m-1}$ satisfying (4.6) then (4.5) can be
  rewritten as
  $
d(I-{\Cal P}_\varepsilon,U)=d(f,U)-m.
  $ Analogously, if $D=0$ then ${\roman s\roman i\roman g\roman n}(q(\theta_j))(2D-1)=1$ is equivalent to $(-1)^{(m-1)/2}{\roman s\roman i\roman g\roman n}(\cos(j\pi))=-1$ that in combination  with (4.4)
  gives ${\roman i\roman n\roman d}(\theta_j,M)=-1$ allowing to rewrite
  (4.5) in the form $
d(I-{\Cal P}_\varepsilon,U)=d(f,U)+m.
  $

\qed

At the end of the paper we note that system (1.2) should exhibit
very complex behavior in order that assumption (4.1) be not
satisfied with any $v_{n-1}\in{\bold R}^{n-1}.$ Particularly,
(4.1) holds true for the prototypic unperturbed system (1.2)

\centerline{$
  \dot x_1=x_2-x_1(x_1^2+x_2^2-1),
$}

\noindent $$
  \dot x_2=-x_1-x_2(x_1^2-x_2^2-1),
  \tag 4.8
$$
\centerline{$
  \dot x_3=-x_3 $}
possessing the nondegenerate $2\pi$-periodic cycle $x_0(t)={\sin
t\choose \cos t}$ and $U=B_1(0)=\{v\in{\bold R}^3:\|v\|<1\}.$
Indeed, it can be easily checked that $
  \Phi(t)=\left({\sin t \choose 0},{\cos t \choose 0},{0\choose 1}\right)^*
  ,$ ${\roman e}^{\Lambda
  t}=\left({{\roman e}^{2t} \atop 0}\ {0\atop {\roman e}^t}\right)$ and  $Y_{n-1}(t)=\Phi(t){\roman
   e}^{-\Lambda t}
$ in this case. Thus, taking $v_{n-1}={1\choose 0}$ we have
$$
  Y_{n-1}(t)\left({\roman e}^{\Lambda
     T}\right)^*\left(\left({\roman
     e}^{\Lambda T}\right)^*-I\right)^{-1}\left({\roman e}^{\Lambda
     t}\right)^*v_{n-1}=\frac{{\roman e}^{2t}}{{\roman
     e}^{2t}-1}\left(\sin t,\cos t,0\right)^*.
$$
This last vector centered at $x_0(t)$ is perpendicular to
$\partial U$ for any $t\in[0,2\pi].$

\enddocument